\documentclass[11pt]{amsart}
\usepackage{amssymb}
\usepackage{graphicx}                            
\usepackage{hyperref}
\usepackage{setspace}                            

\newtheorem{lemma}{Lemma}[section]
\newtheorem{theorem}{Theorem}[section]

\newcommand{\adfmod}[1]{~(\mathrm{mod}~#1)}
\newcommand{\adfhide}[1]{}

\begin{document}
\title{The design spectrum of the Shrikhande graph}
\author{A. D. Forbes}
\address{LSBU Business School,
London South Bank University,
103 Borough Road,
London SE1 0AA, UK.}
\email{anthony.d.forbes@gmail.com}
\author{C. G. Rutherford}
\address{LSBU Business School,
London South Bank University,
103 Borough Road,
London SE1 0AA, UK.}
\email{c.g.rutherford@lsbu.ac.uk}
\date{\today\ version 2.2}

\begin{abstract}
The design spectrum of a simple graph $G$ is the set of positive integers $n$ such that there exists an edgewise
decomposition of the complete graph $K_n$ into $n(n - 1)/(2 |E(G)|)$ copies of $G$.
The purpose of this short paper is to prove that the Shrikhande graph and the line graph of $K_{4,4}$ have the
design spectrum $\{96t + 1: t = 1, 2, \dots\}$.
\end{abstract}

\maketitle


\section{Introduction}\label{sec:Introduction}
If $F$ and $G$ are simple graphs,
an {\em edgewise decomposition} of $F$ into $G$, which we also refer to as
a {\em $G$-decomposition of $F$}, is a partition $\mathcal{E}$ of the edges of $F$
such that each $E \in \mathcal{E}$ is the edge set of a graph isomorphic to $G$.
If $F$ is the complete graph $K_n$, we usually refer to the decomposition as a {\em $G$-design} of order $n$.
The {\em design spectrum} of $G$ is the set of positive integers $n$ for which a $G$-design of order $n$ exists.
If $G$ is $d$-regular, the necessary conditions for the existence of a $G$-design are
\begin{align}
n &\ge |V(G)| \text{~or~} n = 1, \nonumber \\
n(n - 1) &\equiv 0 \adfmod{2|E(G)|}, \label{eqn:necessary-conditions} \\
   n - 1 &\equiv 0 \adfmod{d}. \nonumber
\end{align}

Given a $d$-regular graph $G$, by a theorem of Wilson, \cite{Wilson1976}, the conditions (\ref{eqn:necessary-conditions})
are sufficient for all sufficiently large $n$, and hence
the determination of $G$'s design spectrum is actually a finite problem.
However, it is usually impossible to resolve all of the cases not covered by `sufficiently large'
whenever $d$ or the chromatic number is large.
Nevertheless, design spectra have been computed for many graphs, including some infinite classes, \cite{AdamsBryantBuchanan2008}.
In particular, the design spectrum has been resolved for the Petersen graph, \cite{AdamsBryant1996}.

\begin{figure}
\caption{The Shrikhande graph}
\begin{center}
\includegraphics[width=0.8\textwidth,trim=0mm 0mm 0mm 0mm, clip]{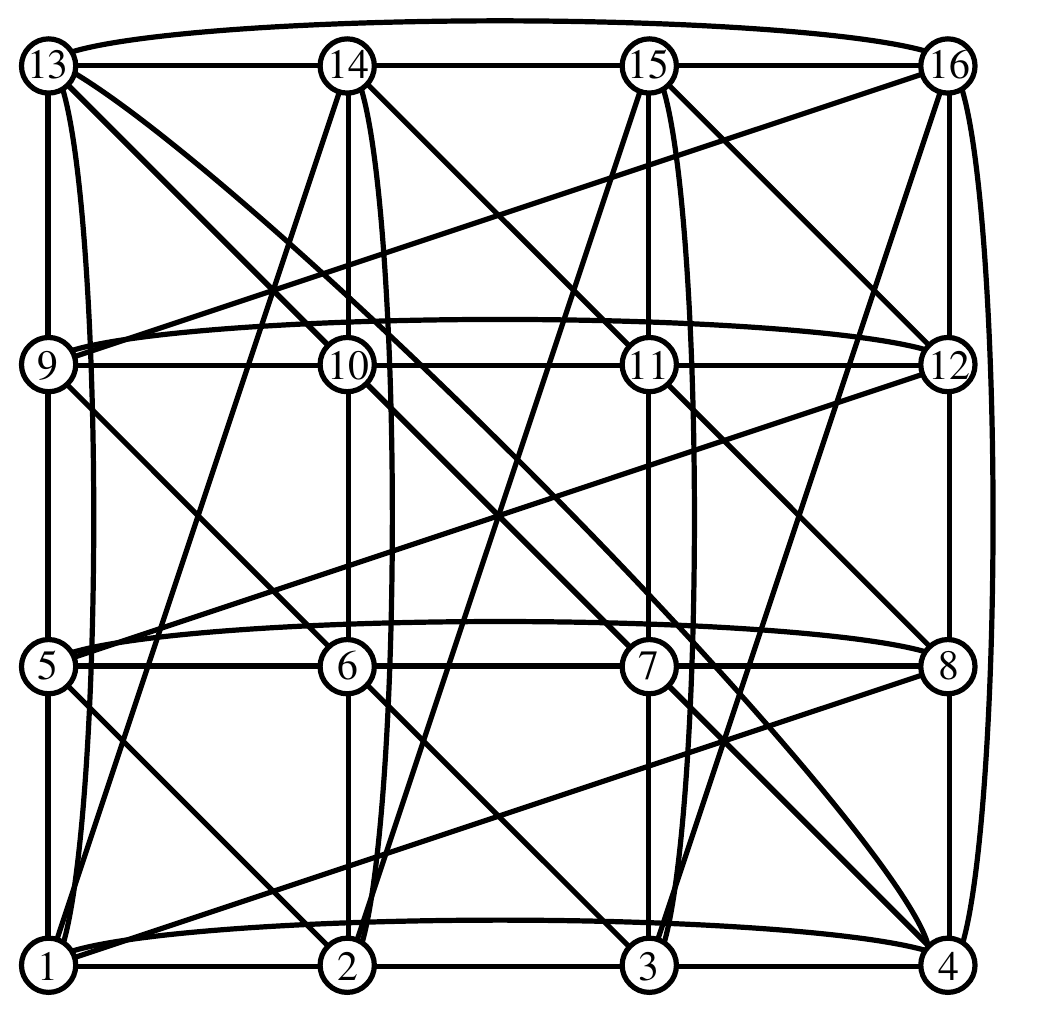} 
\end{center}
\label{fig:Shrikhande-graph}
\end{figure}

In our paper we address the design spectrum problem for the Shrikhande graph, illustrated in Figure~\ref{fig:Shrikhande-graph}, as well as the line graph of $K_{4,4}$.
Our objective is to prove in Theorems~\ref{thm:Shrikhande-graph-design-spectrum} and \ref{thm:LineK44-graph-design-spectrum} that in the conditions (\ref{eqn:necessary-conditions}) are sufficient in each case.
Incidentally, we mention that there exists a decomposition into 5 Shrikhande graphs of $2K_{16}$;
see \cite{Bryant2013}.

The Shrikhande graph is strongly regular with parameters srg(16,6,2,2) and eigenvalues $6^1 2^6 (-2)^9$,
it is 4-chromatic,
and its complement is 6-chromatic.
It is named after Sharad--Chandra S. Shrikhande, who showed that the graph is identified by these properties, \cite{Shrikhande1959}.
For a biographical account of Shrikhande, various constructions of the Shrikhande graph, and
a textbook account of the relevant areas of discrete mathematics, see \cite{CameronVijayakumarLakshmaanan2025}.

\section{The design spectrum of the Shrikhande graph}\label{sec:Shrikhande}
For our proof of Theorem~\ref{thm:Shrikhande-graph-design-spectrum},
we employ a technique of design theory known as Wilson's fundamental construction, \cite{Wilson1972}.
The method uses group divisible designs to build large graph decompositions from small ones.

For the purpose of this paper, a {\em group divisible design}, $k$-GDD, of type $g^{u}$
is an ordered triple ($V, \mathcal{G}, \mathcal{B}$)
where
\begin{enumerate}
\item[(i)]$V$ is a set of $g u$ {\em points},
\item[(ii)]$\mathcal{G}$ is a partition of $V$ into $u$ subsets of size $g$, called \textit{groups}, and
\item[(iii)]$\mathcal{B}$ is a collection of $k$-subsets $V$, called \textit{blocks},
    which has the property that each pair of points from distinct groups occurs in precisely one block but a pair of distinct points from the same group does not occur in any block.
\end{enumerate}

Our first lemma asserts the existence of the group divisible designs that we require
for our constructions.
%
\begin{lemma}\label{lem:GDD-existence}~
There exist $4$-$\mathrm{GDD}$s of type $24^{t}$ for $t \ge 4$.
\end{lemma}
%
\begin{proof}
See \cite{BrouwerSchrijverHanani1977}.
\end{proof}

Next, we give the direct constructions of decompositions into the Shrikhande graph
that we require for our proof of Theorem~\ref{thm:Shrikhande-graph-design-spectrum}.
The sets of labelled graphs that form the decompositions are created from base blocks.
A {\em base block} is an ordered 16-tuple
$(\ell_1, \ell_2, \dots, \ell_{16})$
where, for $i \in \{1, 2, \dots, 16\}$,
label $\ell_i$ is attached to vertex $i$ of the Shrikhande graph as defined by the edge set
\begin{align*}
\{&\{1,2\}, \{1,4\}, \{1,5\}, \{1,8\}, \{1,13\}, \{1,14\}, \{2,3\}, \{2,5\}, \\
  &\{2,6\}, \{2,14\}, \{2,15\}, \{3,4\}, \{3,6\}, \{3,7\}, \{3,15\}, \{3,16\}, \\
  &\{4,7\}, \{4,8\}, \{4,13\}, \{4,16\}, \{5,6\}, \{5,8\}, \{5,9\}, \{5,12\}, \\
  &\{6,7\}, \{6,9\}, \{6,10\}, \{7,8\}, \{7,10\}, \{7,11\}, \{8,11\}, \{8,12\}, \\
  &\{9,10\}, \{9,12\}, \{9,13\}, \{9,16\}, \{10,11\}, \{10,13\}, \{10,14\}, \{11,12\}, \\
  &\{11,14\}, \{11,15\}, \{12,15\}, \{12,16\}, \{13,14\}, \{13,16\}, \{14,15\}, \{15,16\}\}.
\end{align*}
The technique is explained more fully in \cite{ForbesRutherford2025}.
%
%
%
%
%
%
%
\begin{lemma} \label{lem-Design-4-4-4-4}
There exists an edgewise decomposition of the complete $4$-partite graph $K_{4,4,4,4}$ into $2$ copies of the Shrikhande graph.
\end{lemma}
%
\begin{proof}
The point set is $\mathbb{Z}_{16}$ partitioned by residue class modulo 4.
The decomposition consists of two base blocks:


$(0,1,2,5,3,4,7,6,10,9,8,13,11,14,15,12)$,

$(0,2,8,10,9,3,5,15,12,14,4,6,7,13,11,1)$.
%
\end{proof}
%
\begin{lemma} \label{lem-Design-97-et-al}
There exist Shrikhande-graph-designs of orders
$97$, $193$ and $289$.
\end{lemma}
%
\begin{proof}
The graphs for the design of order $n$ are developed from a single base block by
$x \mapsto \omega^e x + d$, $0 \le  e < (n - 1)/96$, $0 \le d < n$,
where $\omega$ is a specified parameter.
For order $n \in \{97, 193\}$, the arithmetic is performed in the field $\mathbb{Z}_{n}$.
For order 289, we use $\mathrm{GF}(17^2)$,
where element $az + b$, $a, b \in \mathbb{Z}_{17}$, is represented by the number $17a + b$.
The polynomial for multiplication in $\mathrm{GF}(17^2)$ is $z^2 + 3z + 1$.

{Design order 97}, $\omega = 1$:


$(0,4,6,62,1,11,19,45,69,80,59,78,32,74,28,44)$;


{Design order 193}, $\omega = 81$:


$(0,19,164,27,51,175,66,138,74,20,70,94,108,77,41,134)$;


{Design order 289}, $\omega = 139$:


$(0,136,232,11,176,180,89,159,288,257,90,42,45,260,37,19)$.
%
%
\end{proof}

\begin{theorem}\label{thm:Shrikhande-graph-design-spectrum}
There exists a Shrikhande-graph-design of order $n \ge 1$ if and only if $n \equiv 1 \adfmod{96}$.
\end{theorem}
%
\begin{proof}
A straightforward computation confirms that the necessary conditions (\ref{eqn:necessary-conditions})
for the existence of a Shrikhande-graph-design of order $n$ simplify to
$n = 96t + 1$, $t = 1, 2, \dots$.

Take a 4-GDD of type $24^{t}$ from Lemma~\ref{lem:GDD-existence},
inflate its points by a factor of 4 and replace its blocks by decompositions into Shrikhande graphs of $K_{4,4,4,4}$,
which exist by Lemma~\ref{lem-Design-4-4-4-4}.
Add a new point and overlay each group plus the new point with a Shrikhande-graph-design of order 97 from Lemma~\ref{lem-Design-97-et-al}.
The result is a Shrikhande-graph-design of order $96t + 1$ for $t \ge 4$.

The values of $n$ not accounted for by the construction are 1, 97, 193 and 289.
For $n = 1$ the design is the empty set.
The other designs are provided by Lemma~\ref{lem-Design-97-et-al}.
\end{proof}

\section{The design spectrum of $L(K_{4,4})$}\label{sec:LineK44}%
The line graph of $K_{4,4}$, denoted by $L(K_{4,4})$, has the same eigenvalue spectrum as the Shrikhande graph.
Here we show that it has the same design spectrum.
%
\begin{theorem}\label{thm:LineK44-graph-design-spectrum}
There exists an $L(K_{4,4})$-design of order $n \ge 1$ if and only if $n \equiv 1 \adfmod{96}$.
\end{theorem}
%
\begin{proof}
The details are as set out in Section~\ref{sec:Shrikhande} for the Shrikhande graph.
It suffices only to specify the edge set of $L(K_{4,4})$ that we used for our computations,
\{\{1,6\}, \{1,7\}, \{1,9\}, \{1,11\}, \{1,15\}, \{1,16\}, \{2,7\}, \{2,8\}, \{2,10\}, \{2,12\}, \{2,13\}, \{2,16\}, \{3,5\}, \{3,8\}, \{3,9\}, \{3,11\}, \{3,13\}, \{3,14\}, \{4,5\}, \{4,6\}, \{4,10\}, \{4,12\}, \{4,14\}, \{4,15\}, \{5,10\}, \{5,11\}, \{5,13\}, \{5,15\}, \{6,11\}, \{6,12\}, \{6,14\}, \{6,16\}, \{7,9\}, \{7,12\}, \{7,13\}, \{7,15\}, \{8,9\}, \{8,10\}, \{8,14\}, \{8,16\}, \{9,14\}, \{9,15\}, \{10,15\}, \{10,16\}, \{11,13\}, \{11,16\}, \{12,13\}, \{12,14\}\}%
, as well as the base blocks for the four decompositions of Lemmas~\ref{lem-Design-4-4-4-4} and \ref{lem-Design-97-et-al}:

{Decomposition of $K_{4,4,4,4}$}:

$(0,1,2,3,4,6,10,8,7,14,5,12,11,9,13,15)$,

$(0,4,5,13,8,2,9,7,14,6,3,15,10,12,11,1)$;
%

{Design order 97}, $\omega = 1$:


$(0,49,34,94,1,2,6,76,68,58,47,80,59,83,38,75)$;


{Design order 193}, $\omega = 81$:


$(0,78,184,67,34,169,130,177,64,103,137,108,30,7,84,65)$;


{Design order 289}, $\omega = 139$:


$(0,119,123,283,273,249,231,69,2,234,93,67,55,64,20,256)$;
%
%
\end{proof}

Finally we have the following.
%
\begin{theorem}\label{thm:Shrikhande+LineK44=K4444}
There exists a decomposition of the complete multipartite graph $K_{4,4,4,4}$
into a Shrikhande graph and a line graph $L(K_{4,4})$.
\end{theorem}
%
\begin{proof}
Combine the edges of the Shrikhande graph as depicted in Figure~\ref{fig:Shrikhande-graph}
with the edges of the line graph $L(K_{4,4})$ as specified in the proof of Theorem~\ref{thm:LineK44-graph-design-spectrum}.
The result is a graph that is isomorphic to $K_{4,4,4,4}$.
\end{proof}


\section*{ORCID}

\noindent A. D. Forbes     \url{https://orcid.org/0000-0003-3805-7056} \\
C. G. Rutherford           \url{https://orcid.org/0000-0003-1924-207X}


\adfhide{
\section{Declarations}

\subsection{Funding}

The authors have no relevant financial or non-financial interests to disclose.

\subsection{Competing Interests}

The authors declare that no funds, grants, or other support were received during the preparation of this manuscript.

\subsection{Data Availability}

We do not analyse or generate any datasets.
}

\end{document}